\documentclass[reqno,centertags, 12pt]{amsart}
\usepackage{amsmath,amsthm,amscd,amssymb}
\usepackage{latexsym}
\newcommand{\bbU}{{\mathbb{U}}}

\newcommand{\bbD}{{\mathbb{D}}}

\newcommand{\bbZ}{{\mathbb{Z}}}
\newcommand{\bbC}{{\mathbb{C}}}
\newcommand{\bbH}{{\mathbb{H}}}

\newcommand{\calP}{{\mathcal P}}

\newcommand{\dott}{\,\cdot\,}

\newcommand{\lb}{\label}
\newcommand{\f}{\frac}

\newcommand{\ol}{\overline}
\newcommand{\ti}{\tilde  }

\newcommand{\s}{\text{\rm{s}}}

\newcommand{\bi}{\bibitem}

\newcommand{\beq}{\begin{equation}}
\newcommand{\eeq}{\end{equation}}
\newcommand{\ba}{\begin{align}}
\newcommand{\ea}{\end{align}}
\newcommand{\veps}{\varepsilon}
\newcounter{smalllist}
\newenvironment{SL}{\begin{list}{{\rm\roman{smalllist})}}{%
\setlength{\topsep}{0mm}\setlength{\parsep}{0mm}\setlength{\itemsep}{0mm}%
\setlength{\labelwidth}{2em}\setlength{\leftmargin}{2em}\usecounter{smalllist}%
}}{\end{list}}

\newcommand{\bigtimes}{\mathop{\mathchoice%
{\smash{\vcenter{\hbox{\LARGE$\times$}}}\vphantom{\prod}}%
{\smash{\vcenter{\hbox{\Large$\times$}}}\vphantom{\prod}}%
{\times}%
{\times}%
}\displaylimits}

\DeclareMathOperator{\Real}{Re}

\allowdisplaybreaks
\numberwithin{equation}{section}
\newtheorem{theorem}{Theorem}[section]
\newtheorem{proposition}[theorem]{Proposition}
\newtheorem{lemma}[theorem]{Lemma}
\newtheorem{corollary}[theorem]{Corollary}
\theoremstyle{definition}

\theoremstyle{remark}

\newcommand{\abs}[1]{\lvert#1\rvert}

\begin{document}
\title{The Sharp Form of the Strong Szeg\H{o} Theorem}
\author{Barry Simon}

\thanks{$^1$ Mathematics 253-37, California Institute of Technology, Pasadena, CA 91125. 
E-mail: bsimon@caltech.edu}
\thanks{$^2$ Supported in part by NSF grant DMS-0140592} 
\thanks{{\it To appear in Proc. Conf. on Geometry and Spectral Theory}}

\dedicatory{In memoriam, Robert Brooks {\rm{(}}1952--2002{\rm{)}} }

\date{January 21, 2004}

\begin{abstract} Let $f$ be a function on the unit circle and $D_n(f)$ be the determinant 
of the $(n+1)\times (n+1)$ matrix with elements $\{c_{j-i}\}_{0\leq i,j\leq n}$ where $c_m  
=\hat f_m\equiv \int e^{-im\theta} f(\theta) \f{d\theta}{2\pi}$. The sharp form of the 
strong Szeg\H{o} theorem says that for any real-valued $L$ on the unit circle with 
$L,e^L$ in $L^1 (\f{d\theta}{2\pi})$, we have 
\[
\lim_{n\to\infty}\, D_n(e^L) e^{-(n+1)\hat L_0} = \exp \biggl( \, \sum_{k=1}^\infty 
k\abs{\hat L_k}^2\biggr)
\]
where the right side may be finite or infinite. We focus on two issues here: a new proof 
when $e^{i\theta}\to L(\theta)$ is analytic and known simple arguments that go from 
the analytic case to the general case. We add background material to make this article 
self-contained. 
\end{abstract}

\maketitle

\section{Introduction} \lb{s1}

Let $\{c_m\}_{m=-\infty}^\infty$ be a two-sided sequence of complex numbers. A {\it Toeplitz 
matrix} is a finite matrix constant along diagonals: 
\begin{equation} \lb{1.1} 
T_{n+1} = \begin{pmatrix} 
c_0 & c_1 & c_2 & \dots & c_n \\
c_{-1} & c_0 & c_1 & \dots & c_{n-1} \\
\vdots  & {} & {} & {} & \vdots \\
c_{-n} & c_{-n+1} & {}  & \dots & c_0 
\end{pmatrix} 
\end{equation} 
It turns out that the natural way to label $T$ is in terms of the Fourier transform of 
$c$, that is, 
\begin{equation} \lb{1.2} 
f(\theta) =\sum_{m=-\infty}^\infty c_m e^{im\theta}
\end{equation}
on $\partial\bbD$ ($\bbD=\{z\mid \abs{z}<1\}$, $\partial\bbD =\{z\mid\abs{z}=1\}$). $f$ 
is called the {\it symbol} of the Toeplitz matrix. 

One can define a symbol as a distribution so long as $\abs{c_m}$ is polynomially bounded 
in $m$, but we will discuss the case where there is a signed measure, $d\mu$, so that 
\begin{equation} \lb{1.3} 
c_n =\int e^{-in\theta}\, d\mu(\theta) \equiv \hat\mu_n 
\end{equation}
As usual, for $f\in L^1 (\partial\bbD, \f{d\theta}{2\pi})$, we define $\hat f_n$ to be the 
Fourier coefficients of the measure $f\f{d\theta}{2\pi}$. We will most often discuss the 
case where $d\mu$ is absolutely continuous, that is, $d\mu = w(\theta)\f{d\theta}{2\pi}$ 
and where $w\geq 0$ or even that $w=e^L$. 

$D_n(d\mu)$ is the determinant of $T_{n+1}$. The strong Szeg\H{o} theorem says that if $L, 
e^L\in L^1$ with $L$ real, then 
\begin{equation} \lb{1.4} 
\log D_n \biggl(e^L \,\f{d\theta}{2\pi}\biggr) \sim (n+1) \hat L_0 + \sum_{k=1}^\infty\, 
\abs{k} \, \abs{\hat L_k}^2
\end{equation}

There are a number of remarkable aspects of \eqref{1.4}. The first term was found in 1915 
and the second in 1952. Despite the 37-year break, they were both found by Szeg\H{o} --- 
the twenty-year old in 1915 \cite{Sz15} and the 57-year old in 1952 \cite{Sz52}! You 
might wonder about whether \eqref{1.4} is the leading term in a systematic $1/n$ series. 
In fact, if $L$ is real-valued and if $e^{i\theta}\mapsto L(\theta)$ is analytic in 
the neighborhood of $\partial\bbD$, then the error in \eqref{1.4} is $O(e^{-Bn})$ --- 
there are no more terms in the series (this follows from \eqref{2.21new} and \eqref{5.16} 
below). Lest you be shocked by this, we note that for many models in statistical mechanics, 
the free energy has a volume term, a surface term, and then, if the interaction is 
short-range, exponentially small errors. 

A second remarkable aspect is the subtlety. Why should $\log w$ enter at all, and 
then in both linear and quadratic terms? There is a fascination with this subject 
among mathematicians who have extended the result both in the context of function 
algebras \cite{Werm,HoffActa,Lum} and in the context of pseudodifferential operators 
on manifolds \cite{W79,GOk2} (see \cite{OPUC} for literally dozens of papers on 
each aspect). 

A third aspect is that there are a remarkable number of applications of this result. 
Szeg\H{o} returned to find the second term because of a question raised by Onsager 
who ran into Toeplitz determinants in his work on the Ising model (see \cite{McCoy,BoOns} 
for a discussion of this). They enter in the study of some Coulomb gases 
\cite{Len64,Len72,FH69} and in electrical engineering applications \cite{Kai74,BakCon}. 
And they have had a surge of interest recently because of their role in random matrix 
theory \cite{MehBk}. 

When Szeg\H{o} \cite{Sz52} proved \eqref{1.4}, he assumed $L$ was $C^{1+\veps}$. There 
were many papers on this subject which improved this incrementally until Ibragimov  
\cite{Ib}, fifteen years later, proved the following sharp form: 

\begin{theorem}[\cite{Ib,GoIb}]\lb{T1.1} Let $L$ be a real-valued function on $\partial\bbD$ 
so that $L,e^L\in L^1 (\partial\bbD, \f{d\theta}{2\pi})$. Then 
\begin{equation} \lb{1.5} 
\lim_{n\to\infty}\, D_n \biggl( e^L \, \f{d\theta}{2\pi}\biggr) e^{-(n+1)\hat L_0} = 
\exp \biggl(\, \sum_{k=1}^\infty\, \abs{k}\, \abs{\hat L_k}^2\biggr)
\end{equation}
\end{theorem} 

Thus, \eqref{1.4} holds whenever the right side makes sense, that is, $L$ in $H^{1/2}$, 
the Sobolev space of order $\f12$. This should be supplemented with a result of 
Golinskii-Ibragimov \cite{GoIb}: 

\begin{theorem}\lb{T1.2} If $d\mu=e^{L}\f{d\theta}{2\pi} + d\mu_\s$ with $d\mu_\s$ singular 
is a positive measure on $\partial\bbD$ and $\lim_{n\to\infty} D_n (d\mu) e^{-(n+1)\hat L_0} 
<\infty$, then $d\mu_\s =0$. 
\end{theorem} 

The combination of these two theorems has a spectral theory consequence that 
links it up to the theme of the conference and to Bob Brooks' interests. As we will 
see in Section~\ref{s2}, probability measures on $\partial\bbD$ have associated 
parameters $\{\alpha_n (d\mu)\}_{n=0}^\infty$ called Verblunsky coefficients. 
It can be shown using Theorems~\ref{T1.1}, \ref{T1.2}, and \ref{T2.4} that 

\begin{theorem}\lb{T1.3} Let $d\mu$ be a probability measure on $\partial\bbD$ and 
$\{\alpha_n(d\mu)\}_{n=0}^\infty$ its Verblunsky coefficients. Then the following are 
equivalent: 
\begin{SL} 
\item[{\rm{(i)}}] $\sum_{n=0}^\infty n\abs{\alpha_n}^2 <\infty$ 
\item[{\rm{(ii)}}] $d\mu_\s =0$ and $d\mu =e^L \f{d\theta}{2\pi}$ where $\sum_{k=1}^\infty 
k \abs{\hat L_k}^2 <\infty$. 
\end{SL} 
\end{theorem} 

This is one of those gems of spectral theory that give necessary and sufficient conditions 
relating properties of a measure and its inverse spectral parameters. 

We have two main themes in this article. First, we wish to note that despite it taking 
fifteen years to go from $L\in C^{1+\veps}$ to $L\in H^{1/2}$, there is an elegant 
and simple argument to do this jump. This combines arguments of Golinskii-Ibragimov 
\cite{GoIb} and Johansson \cite{Jo88} whose proofs have ``easy" halves that handle 
opposite sides of the extension. It does not appear to be widely appreciated that 
their arguments can be combined in this way. In fact, these results reduce \eqref{1.4} 
to the case where the $\hat L_k$ decay exponentially, that is, $e^{i\theta} \mapsto 
L(\theta)$ is real analytic on $\partial\bbD$. 

Second, we have a new proof of \eqref{1.4} in this real analytic case that is perhaps 
less mysterious than the elaborate calculation in Szeg\H{o} \cite{Sz52}. From our 
point of view, the two terms in \eqref{1.4} come from two terms in the Christoffel-Darboux 
formula for $z=e^{i\theta}$. 

While the arguments of \cite{GoIb,Jo88} are simple, they depend on considerable 
general machinery relating Toeplitz determinants to orthogonal polynomials on the 
unit circle (OPUC) on the one hand, and to the statistical mechanics of Coulomb 
gases on the other (both themes go back to Szeg\H{o}'s early work: \cite{Sz15} 
has the Coulomb gas representation and the main point of \cite{Sz20,Sz21} is to 
discuss the connection to OPUC), so this article, in attempting to be self-contained, 
provides this background.  

Sections~\ref{s2} and \ref{s3} discuss the basics of OPUC. In Section~\ref{s4}, we get the 
leading term in \eqref{1.4}, not only for its own sake, but to define in Section~\ref{s5} 
the Szeg\H{o} function which will play a critical role in many aspects of the remainder. 
These preliminaries allow us to present the Golinskii-Ibragimov half of the extension in 
Section~\ref{s6}. After proving the Coloumb gas representation in Section~\ref{s7}, we 
can prove the Johansson half of the extension in Section~\ref{s8}. The final three 
sections provide the proof of \eqref{1.4} in the analytic case: Section~\ref{s9} has a 
preliminary proving the Christoffel-Darboux formula, and the last two sections finish 
the proof. 

I have written a comprehensive book on OPUC \cite{OPUC} and everything in this paper 
appears there, but it seemed sensible, given the fact that the material is spread 
through a long book, to pull out exactly what is needed to prove Ibragimov's theorem. 

While \eqref{1.4} is sharp in one sense, it is not the end of the story by any means. 
First, there is a simple argument of Johansson \cite{Jo88} that drops the requirement 
of reality from $L$: if $e^L,L\in L^1$ and $L\in H^{1/2}$, then an extension \eqref{1.4} 
holds in the sense that 
\[
e^{-(n+1)\hat L_0} D_n \biggl( e^L \, \f{d\theta}{2\pi}\biggr) \to \exp 
\biggl(\, \sum_{k=1}^\infty\, \abs{k} \hat L_k \hat L_{-k}\biggr)
\]
There are also subtleties in extending \eqref{1.4} to allow matrix-valued symbols, 
to allow complex $w$'s with nonzero winding number, and to determine the leading 
behavior when $L\notin L^1$. The reader can consult \cite{OPUC} for references 
on all these issues. 

\smallskip
Over the course of studying asymptotics of Toeplitz determinants, I have learned 
a lot in discussions with Percy Deift, Rowan Killip, and Irina Nenciu, and I would 
like to thank them for their insights. 

Bob Brooks was a substantial mathematician and wonderful person. We lost him 
too soon. I'm glad to dedicate this article to his memory.

%%%%%%%%%%%%%%%%%%%%%%%%%%%%%%%%%%%%%%%%%%%%%%%%%%%%%%%%%%%%%%%%
\section{Verblunsky Coefficients and Toeplitz Matrices} \lb{s2}
%%%%%%%%%%%%%%%%%%%%%%%%%%%%%%%%%%%%%%%%%%%%%%%%%%%%%%%%%%%%%%%%

If $c_n$ are the moments of a measure, $\mu$, it is natural to form the monic orthogonal 
polynomials, $\Phi_n (z;d\mu)$, defined by
\begin{gather} 
\Phi_n(z) =z^n + \text{ lower order}  \lb{2.1} \\
\langle z^j, \Phi_n\rangle_\mu =0 \qquad \text{if } j=0,1,\dots, n-1 \lb{2.2}
\end{gather} 
where 
\begin{equation} \lb{2.3} 
\langle f,g\rangle = \int\, \ol{f(e^{i\theta})}\, g(e^{i\theta})\, d\mu(\theta)
\end{equation}
is the $L^2 (\partial\bbD, d\mu)$ inner product. In order to form $\Phi_n$ for all 
$n$, we need the measures $d\mu$ to be nontrivial, that is, not supported on a 
finite set of points. 

The matrix elements $c_{k-\ell}=\int e^{i(\ell-k)\theta}\ d\mu =\langle z^k, z^\ell 
\rangle_\mu$, so $D_n (d\mu)$ is a Gram determinant. Such determinants allow change 
of basis, that is, if $P_k(z)=z^k + \text{ lower order}$, $\det (\langle P_k,P_\ell
\rangle )_{0\leq k,\ell\leq n} =D_n (d\mu)$. We can take $P_k =\Phi_k$, in which 
case, $\langle P_k,P_\ell\rangle$ is a  diagonal matrix (!), and so, 

\begin{theorem}\lb{T2.1} 
\begin{equation} \lb{2.3a} 
D_n (d\mu) =\prod_{j=0}^n \, \|\Phi_j\|_{d\mu}^2
\end{equation} 
\end{theorem} 

$\Phi_j$ is the orthogonal projection in $L^2$ of $z^j$ onto the orthogonal complement 
of  $[1,\dots, z^{j-1}]$, the span  of $\{1, \dots, z^{j-1}\}$. Since multiplication 
by $z$ is unitary, $z\Phi_j$ is the projection of $z^{j+1}$ onto $[z,\dots, z^j]^\perp$ 
while $\Phi_{j+1}$ is the projection of $z^{j+1}$ onto $[1,\dots, z^j]^\perp$, so  
\begin{equation} \lb{2.3b} 
\|\Phi_{j+1}\| \leq \|z\Phi_j\| = \|\Phi_j\|
\end{equation}
Thus, $\|\Phi_j\|$ is decreasing in $j$. It follows that 

\begin{theorem} \lb{T2.2} 
\begin{alignat}{2} 
&\text{\rm{(a)}} \qquad && \lim\, D_n (d\mu)^{1/n+1} = \lim \, 
\f{D_{n+1}(d\mu)}{D_n (d\mu)} = \lim_{n\to\infty}\, \|\Phi_n\|^2 \lb{2.4} \\
&\text{\rm{(b)}} \qquad && \f{D_{n+1}}{D_n} \leq \f{D_n}{D_{n-1}} \lb{2.5} 
\end{alignat} 
\end{theorem} 

These ideas all go back to Szeg\H{o} \cite{Sz20,Sz21}. The next idea, which is a 
set of recursion relations of the $\Phi_n$, was first written down by Szeg\H{o} 
\cite{Szb}, but the basic parameters occurred  in a related context in Verblunsky 
\cite{V35,V36}. To state them, we need to define the reversed polynomials 
\begin{equation} \lb{2.6} 
\Phi_n^*(z) = z^n \, \ol{\Phi_n (1/\bar z)}
\end{equation} 

\begin{theorem}\lb{T2.3} For any nontrivial measure, $d\mu$, there exists a 
sequence of numbers $\{\alpha_n\}_{n=0}^\infty$ so that 
\begin{equation} \lb{2.7} 
\Phi_{n+1}(z) =z\Phi_n(z) -\bar\alpha_n \Phi_n^*(z)
\end{equation}
Moreover, $\alpha_n\in\bbD$ and 
\begin{equation} \lb{2.8} 
\|\Phi_{n+1}\|^2 = (1-\abs{\alpha_n}^2) \|\Phi_n\|^2
\end{equation} 
and 
\begin{equation} \lb{2.8a} 
\|\Phi_n\|^2 = c_0 (d\mu) \prod_{j=0}^{n-1} \, (1-\abs{\alpha_j}^2)
\end{equation} 
\end{theorem} 

{\it Remarks.} 1. In the next section, we will see that $\mu\mapsto 
\{\alpha_n\}_{n=0}^\infty$ is one-one for $\mu$'s which are normalized. 

\smallskip
2. There is a converse going back to Verblunsky \cite{V35} (see \cite{OPUC} for 
many other proofs) that the map from probability measures to $\bigtimes_{n=0}^\infty 
\bbD$ by $\mu\mapsto \{\alpha_n\}_{n=0}^\infty$ is onto. 

\smallskip
3. The $\alpha_n$ are called the {\it Verblunsky coefficients} for $d\mu$.  

\smallskip
4. Applying $Q^*(z)=z^{n+1}\, \ol{Q(1/\bar z)}$ to \eqref{2.7} yields 
\begin{equation} \lb{2.9} 
\Phi_{n+1}^*(z) = \Phi_n^*(z) -\alpha_n z\Phi_n(z)
\end{equation}

\smallskip
5. The proof below is a variant of one of Atkinson \cite{Atk64}. Szeg\H{o}'s proof 
first proves the Christoffel-Darboux formula (see Section~\ref{s9}) and uses that 
to prove \eqref{2.7}. 

\begin{proof} Let $V_ng =z^n \bar g$ on $L^2 (d\mu)$. $V_n$ is anti-unitary, maps 
$\calP_n$, the polynomials of degree $n$, to themselves, and maps $\Phi_n$ to 
$\Phi_n^*$. Since $\Phi_n$ is the unique element (up to constants) of $\calP_n$ 
orthogonal to $\{1,z,\dots, z^{n-1}\}$ and $V_n$ is anti-unitary, $\Phi_n^*$ is 
the unique element of $\calP_n$ orthogonal to $\{V_n 1, \dots, V_n z^{n-1}\} = 
\{z^n, z^{n-1}, \dots, z\}$. 

Now, for $j=1, \dots, n$, 
\[
\langle z^j, z\Phi_n\rangle = \langle z^{j-1}, \Phi_n\rangle =0 
\]
and clearly, $\langle z^j, \Phi_{n+1}\rangle =0$. Thus, 
\[
\langle z^j, \Phi_{n+1} -z\Phi_n\rangle =0 
\]
for $j=1, \dots, n$. Since $\Phi_n$ and $\Phi_{n+1}$ are monic, $\Phi_{n+1}-z\Phi_n 
\in\calP_n$. So, by the first part of the proof, \eqref{2.7} holds for a suitable 
constant $\alpha_n$. Thus 
\begin{equation} \lb{2.10} 
\alpha_n = -\ol{\Phi_{n+1}(0)}
\end{equation}

Since $\Phi_n^* \perp \Phi_{n+1}$, we have 
\begin{align*} 
\|\Phi_n\|^2 &= \|z\Phi_n\|^2 = \|\Phi_{n+1} + \bar\alpha_n \Phi_n^*\|^2 \\
&= \|\Phi_{n+1}\|^2 + \abs{\alpha_n}^2 \|\Phi_n\|^2 
\end{align*} 
which implies \eqref{2.8}. \eqref{2.8} in turn implies $\abs{\alpha_n}<1$. 
\eqref{2.8a} follows by induction and $\|\Phi_0 \|^2 = \|1\|^2 =c_0 (d\mu)$. 
\end{proof} 

This leads to 

\begin{theorem} \lb{T2.4} Suppose $\int d\mu =1$. We have 
\begin{equation} \lb{2.11} 
F(d\mu) =\lim_{n\to\infty}\, \f{D_{n+1}(d\mu)}{D_n (d\mu)} = \prod_{j=0}^\infty \, 
(1 -\abs{\alpha_j}^2)
\end{equation} 
where the product always converges although the limit may be zero. 

If $F(d\mu) >0$, then 
\begin{equation} \lb{2.12x}  
G_n = \f{D_n}{F^{n+1}}
\end{equation}
obeys 
\begin{equation} \lb{2.12a} 
G_{n+1} \geq G_n
\end{equation}
The limit always exists {\rm{(}}but may be infinite{\rm{)}} and is given by 
\begin{equation} \lb{2.12} 
G(d\mu) =\lim_{n\to\infty}\, G_n(d\mu) =\prod_{j=0}^\infty \, 
(1-\abs{\alpha_j}^2)^{-j-1}
\end{equation}
\end{theorem} 

{\it Remarks.} 1. In particular, 
\begin{align} 
F>0 &\Leftrightarrow \sum_{n=0}^\infty \, \abs{\alpha_n}^2 <\infty \lb{2.13} \\
G<\infty &\Leftrightarrow \sum_{n=0}^\infty \, (n+1) \abs{\alpha_n}^2 < \infty \lb{2.14} 
\end{align} 

\smallskip
2. \eqref{2.11} in a sense goes back to Verblunsky \cite{V36}. \eqref{2.12} seems 
to have only been noted by Baxter \cite{Bax} many years later. 

\smallskip
3. If $c_0\neq 1$, $F(d\mu) =c_0 \prod_{j=0}^\infty (1-\abs{\alpha_j}^2)$ while 
$G(d\mu)$ is still given by \eqref{2.12}. Indeed, $G_n (d\mu) = G_n (d\mu/\int d\mu)$.  

\smallskip
4. \eqref{2.5} says $\log D_n$ is concave in $n$. The monotonicity of $G_n$ is a standard 
fact about concave functions with finite asymptotics. 

\begin{proof}  By \eqref{2.3a} and then \eqref{2.8a}, 
\begin{equation} \lb{2.15} 
\f{D_{n+1}}{D_n} = \|\Phi_{n+1}\| = \prod_{j=0}^n \, (1-\abs{\alpha_j}^2) 
\end{equation}
(if $\|\Phi_1\|^2 =c_0 =1$). From this and $\abs{\alpha_j}<1$, \eqref{2.11} is 
immediate. If $F$ is nonzero, 
\[
\f{D_{n+1}}{D_n}\, \f{1}{F} = \prod_{j=n+1}^\infty\, (1 -\abs{\alpha_j}^2)^{-1} 
\]
so that 
\begin{align} 
G_n &= \f{c_0}{F}\, \prod_{k=1}^n \, \biggl[ \f{D_k}{D_{k-1}} \, \f{1}{F} \biggr] \notag \\
&= \prod_{j=0}^\infty\, (1 -\abs{\alpha_j}^2)^{-\min (n,j)-1} \lb{2.21new} 
\end{align} 
from which $G_{n+1}\geq G_n$ and \eqref{2.12} are immediate. 
\end{proof}

%%%%%%%%%%%%%%%%%%%%%%%%%%%%%%%%%%%%%%%%%%%%%%%%%%%%%%
\section{Bernstein-Szeg\H{o} Approximations} \lb{s3} 
%%%%%%%%%%%%%%%%%%%%%%%%%%%%%%%%%%%%%%%%%%%%%%%%%%%%%% 

Given a nontrivial probability measure, $d\mu$, with Verblunsky coefficients 
$\{\alpha_n\}_{n=0}^\infty$, we will identify the measures, $d\mu^{(N)}$, with
\begin{equation} \lb{3.1} 
\alpha_j (d\mu^{(N)}) = \begin{cases} 
\alpha_j (d\mu) & j=0,1,\dots, N-1 \\
0 & j\geq N \end{cases}
\end{equation}
and see $d\mu^{(N)} \to d\mu$ weakly. In many ways, the general proof of the strong 
Szeg\H{o} theorem will play off this approximation and the distinct approximation 
obtained by truncating the Fourier series for $L$ in $e^L \f{d\theta}{2\pi}$. As a 
preliminary, we need 

\begin{theorem} \lb{T3.1} $\Phi_n (z)$ has all its zeros in $\bbD$. $\Phi_n^*(z)$ has 
all its zeros in $\bbC\backslash\bar\bbD$. 
\end{theorem} 

{\it Remark.} This proof is due to Landau \cite{Land}. See \cite{OPUC} for many other 
proofs of this fact. 

\begin{proof} Let $\Phi_n(z_0) =0$. Then $\pi_{n-1} =\Phi_n(z)/(z-z_0)$ is a polynomial 
of degree $n-1$, and so $\langle\pi_{n-1},\Phi_n\rangle =0$. Since $(z-z_0) \pi_{n-1} 
=\Phi_n$, 
\begin{align*} 
\|\pi_{n-1}\|^2 &= \|z\pi_{n-1}\|^2 = \|z_0 \pi_{n-1} + \Phi_n\|^2 \\
&= \abs{z_0}^2 \|\pi_{n-1}\|^2 + \|\Phi_n\|^2 
\end{align*} 
or 
\begin{equation} \lb{3.2} 
(1-\abs{z_0}^2) \|\pi_{n-1}\|^2 = \|\Phi_n\|^2
\end{equation}
from which we conclude $\abs{z_0} <1$, that is, the zeros of $\Phi_n$ lie in 
$\bbD$. Since $\Phi_n^* (z_0) =0$ if and only if $\Phi_n (1/\bar z_0) =0$, the 
zeros of $\Phi_n^*$ lie in $\bbC\backslash\bar\bbD$. 
\end{proof} 

Given $\mu$, define a measure, $\ti\mu^{(N)}$, by 
\begin{equation} \lb{3.3} 
d\tilde\mu^{(N)} = \f{d\theta}{2\pi \abs{\Phi_N (e^{i\theta};d\mu)}^2}
\end{equation}
which can be defined since $\Phi_N (e^{i\theta})\neq 0$ for all $\theta$ by 
Theorem~\ref{T3.1}. We have 

\begin{lemma} \lb{L3.2} For $j\geq 0$, 
\begin{equation} \lb{3.4} 
\Phi_{N+j} (z;d\ti\mu^{(N)}) = z^j \Phi_N (z;d\mu)
\end{equation}
Moreover, 
\begin{equation} \lb{3.4a} 
\alpha_\ell (d\ti\mu^{(N)}) =0
\end{equation}
for $\ell \geq N$. 
\end{lemma} 

\begin{proof} Since 
\begin{equation} \lb{3.5} 
\abs{\Phi_N(e^{i\theta})}^2 = e^{-iN\theta} \Phi_N (e^{i\theta}) \Phi_N^* (e^{i\theta})
\end{equation}
we have for $k\in\bbZ$, $k\leq N-1$ (includes $k<0)$, 
\begin{align*} 
\int_0^{2\pi} e^{-ik\theta}\, \f{\Phi_N (e^{i\theta})}{\abs{\Phi_N (e^{i\theta})}^2}\, 
\f{d\theta}{2\pi} &= \f{1}{2\pi i}\, \oint_{\abs{z}=1} z^{N-k} \, \f{dz}{z\Phi_N^*(z)} \\
&= 0
\end{align*} 
since $N-k-1 \geq 0$ and $\Phi_N^*(z)^{-1}$ is analytic in a neighborhood of $\bar\bbD$ 
by Theorem~\ref{T3.1}. This says 
\begin{equation} \lb{3.6} 
\langle z^\ell, z^j \Phi_N\rangle_{\ti\mu^{(N)}} =0
\end{equation}
for $\ell=0,1,\dots, N+j-1$. Since $z^j \Phi_N$ is monic, \eqref{3.4} holds. 

By \eqref{2.10}, if $\Phi_{k+1}(0)=0$, then $\alpha_k =0$, so \eqref{3.4} implies 
$\Phi_{N+j} (0;d\ti\mu^{(N)})=0$, which in turn implies \eqref{3.4a}. 
\end{proof} 

\begin{theorem}[Geronimus \cite{Ger46}]\lb{T3.3} Let $d\mu,d\nu$ be two nontrivial 
measures on $\partial\bbD$. Suppose that for some fixed $N$\!, 
\begin{equation} \lb{3.7} 
\Phi_N(z;d\mu) =\Phi_N (z;d\nu)
\end{equation}
Then 
\begin{alignat}{2} 
&\Phi_j (z;d\mu) = \Phi_j (z;d\nu) \qquad && j=0,1,\dots, N-1 \lb{3.8} \\
&\alpha_j (d\mu) = \alpha_j (d\nu) \qquad && j=0,1,\dots, N-1 \lb{3.9} \\
&\f{c_j (d\mu)}{c_0 (d\mu)} = \f{c_j(d\nu)}{c_0(d\nu)} \qquad && j=0,1,\dots, N \lb{3.10} \\
&\|\Phi_j\|_{d\mu}^2 = c_0 (d\mu) c_0 (d\nu)^{-1} \|\Phi_j\|_{d\nu}^2 
\qquad && j=0,1,\dots, N \lb{3.11x} 
\end{alignat} 
\end{theorem} 

\begin{proof} \eqref{2.7} and \eqref{2.9} can be written in matrix form 
\begin{equation} \lb{3.11} 
\begin{pmatrix} \Phi_{j+1}(z) \\ \Phi_{j+1}^*(z) \end{pmatrix} = 
\begin{pmatrix} z & -\bar\alpha_j \\ -z\alpha_j & 1 \end{pmatrix} 
\begin{pmatrix} \Phi_j(z) \\ \Phi_j^*(z) \end{pmatrix}
\end{equation}
The $2\times 2$ matrix in \eqref{3.11} has an inverse $z^{-1} \rho_j^{-2} 
\left(\begin{smallmatrix} 1 & \bar\alpha_j \\ \alpha_j z & z\end{smallmatrix}\right)$ 
where 
\begin{equation} \lb{3.12} 
\rho_j = (1-\abs{\alpha_j}^2)^{1/2}
\end{equation}
Thus \eqref{3.12} implies the inverse Szeg\H{o} recursions: 
\begin{align} 
\Phi_j(z) &= \rho_j^{-2} \f{[\Phi_{j+1}(z) + \bar\alpha_j \Phi_{j+1}^*]}{z} \lb{3.13} \\
\Phi_j^*(z) &= \rho_j^{-2} [\Phi_{j+1}^*(z) + \alpha_j \Phi_{j+1}(z)] \lb{3.14} 
\end{align} 
\eqref{3.7} implies 
\[
\alpha_N (d\mu) = \ol{-\Phi_N(0;d\mu)} = \ol{-\Phi_N(0;d\nu)} = \alpha_N(d\nu) 
\]
and thus, by \eqref{3.13} with $j=N-1$, we have \eqref{3.8} for $j=N-1$. By iterating 
this argument, we conclude that \eqref{3.8} and \eqref{3.9} hold. 

We need only prove \eqref{3.10} and \eqref{3.11x}, assuming $c_0 (d\mu) = c_0 (d\nu) 
=1$ since with $d\ti\mu = d\mu/c_0 (d\mu)$, we have $c_n (d\mu) =c_0 (d\mu)c_n (d\ti\mu)$ 
and $\|\Phi_j\|_{d\mu}^2 =c_0(d\mu) \|\Phi_j\|_{d\ti\mu}^2$. \eqref{3.11x} is immediate 
from \eqref{2.8a}. 

We prove \eqref{3.10} when $c_0(d\mu)=c_0 (d\nu)$ by induction and noting $\langle 1,
\Phi_k\rangle =0$ yields a formula for $c_k$ in terms of the coefficients of $\Phi_k$ 
and $c_0, c_1, \dots, c_{k-1}$. 
\end{proof} 

{\it Remark.} The last paragraph of the proof shows that the $\alpha$'s determine 
the $c$'s and proves that the map of $\mu$ to $\alpha$ is one-one. 

\smallskip

To succinctly state this section's final result, we introduce 
\begin{equation} \lb{3.15} 
\varphi_N (z;d\mu) = \f{\Phi_N (z;d\mu)}{\|\Phi_N\|}
\end{equation} 
the orthonormal polynomials. By \eqref{2.8a}, 
\begin{equation} \lb{3.16} 
\varphi_N (z;d\mu) =\prod_{j=0}^{N-1} \rho_j^{-1} c_0 (d\mu)^{-1/2} \Phi_N (z;d\mu)
\end{equation} 

We define 
\begin{equation} \lb{3.16a} 
\kappa_n = \biggl(\, \prod_{j=0}^{n-1} \rho_j^{-1}\biggr) (c_0 (d\mu))^{-1/2} = 
\|\Phi_n\|^{-1}
\end{equation}
and 
\begin{equation} \lb{3.16b} 
\kappa_\infty = \lim_{n\to\infty}\, \|\Phi_n\|^{-1}
\end{equation}
which exists (but it may be $+\infty$) by $\rho_j <1$. The infinite product in 
\eqref{3.16a}, and so $\kappa_\infty$, is finite if and only if 
\begin{equation} \lb{3.17c} 
\kappa_\infty < \infty \Leftrightarrow \sum_{j=0}^\infty \, \abs{\alpha_j}^2 <\infty 
\end{equation} 

We also note the translation of \eqref{2.7}/\eqref{2.9} from $\Phi$ to $\varphi$: 
\begin{align} 
\varphi_{n+1}(z) &= \rho_n^{-1} (z\varphi_n (z) - \bar\alpha_n \varphi_n^*(z)) \lb{3.17d} \\
\varphi_{n+1}^* (z) &= \rho_n^{-1} (\varphi_n^*(z) - \alpha_n z \varphi_n^*(z)) \lb{3.17e}  
\end{align} 

\begin{theorem}\lb{T3.4} Let $d\mu$ be a nontrivial probability measure on $\partial\bbD$. 
Define 
\begin{equation} \lb{3.17} 
d\mu^{(N)} =\f{d\theta}{2\pi\abs{\varphi_N (e^{i\theta})}^2}
\end{equation}
Then $d\mu^{(N)}$ is a probability measure on $\partial\bbD$ for which \eqref{3.1} 
holds. As $N\to\infty$, $d\mu^{(N)}\to d\mu$ weakly. 
\end{theorem} 

{\it Remark.} We call measures of the form \eqref{3.17} {\it BS measures} and 
$d\mu^{(N)}$ the {\it BS approximation}. 

\begin{proof} Since $d\mu^{(N)}$ is a multiple of $d\ti\mu^{(N)}$, we have the bottom 
half of \eqref{3.1} by \eqref{3.4a}. Since \eqref{3.4} holds for $j=0$, Theorem~\ref{T3.3} 
and \eqref{3.9} imply the top half of \eqref{3.1}.  

Since $\Phi_N =\|\Phi_N\|_\mu \varphi_N$, we clearly have 
\[
\|\Phi_N\|_{\mu^{(N)}}^2 = \|\Phi_N\|_\mu^2
\]
so, by \eqref{3.11}, $c_0 (d\mu^{(N)}) =c_0 (d\mu) =1$, that is, $d\mu_N$ is a probability 
measure. 

By the above and \eqref{3.10}, we have 
\begin{equation} \lb{3.18} 
c_j (d\mu^{(N)}) = c_j (d\mu) \qquad j=0,1,\dots, N
\end{equation}
This and its complex conjugate implies that for any Laurent polynomial, $f$ (polynomial in $z$ 
and $z^{-1}$), 
\begin{equation} \lb{3.19} 
\lim_{N\to\infty}\, \int f(e^{i\theta})\, d\mu^{(N)} = \int f(e^{i\theta})\, d\mu
\end{equation}
since the left side is equal to the right for $N$ large. Since Laurent polynomials are 
dense in $C(\partial\bbD)$, \eqref{3.19} holds for all $f$, that is, we have weak 
convergence. 
\end{proof} 

We note we have proven that 
\begin{equation} \lb{3.20} 
\alpha_j (d\mu^{(N)}) = \begin{cases} \alpha_j & j\leq N-1 \\
0 & j\geq N \end{cases}
\end{equation}

The ideas of this section go back to Geronimus \cite{Ger46} and were rediscovered in 
\cite{ENZG91} and \cite{DGK78}. In particular, the use of inverse recursion to prove 
Geronimus' theorem is taken from \cite{DGK78}.

%%%%%%%%%%%%%%%%%%%%%%%%%%%%%%%%%%%%%%
\section{Szeg\H{o}'s Theorem} \lb{s4} 
%%%%%%%%%%%%%%%%%%%%%%%%%%%%%%%%%%%%%%

In this section, our main goal is to prove 

\begin{theorem}\lb{T4.1} For any $w\in L^1 (\f{d\theta}{2\pi})$, 
\begin{equation} \lb{4.1} 
\lim_{N\to\infty}\, \f{1}{N} \, \log D_N(w) = \int \log w(\theta) \, \f{d\theta}{2\pi} 
\end{equation}
\end{theorem} 

{\it Remarks.} 1. Since $\log(x) \leq x-1$ and $w(\theta)\in L^1$, $\int \max (0,\log
w(x)) \f{d\theta}{2\pi} <\infty$, so the integral on the right side of \eqref{4.1} is 
either convergent or diverges to $-\infty$, in which case \eqref{4.1} says $D_n(w)^{1/n} 
\to 0$. 

\smallskip
2. This was conjectured by P\'olya \cite{Pol} and proven by the twenty-year old 
Szeg\H{o} in 1915 \cite{Sz15}. Our proof here is essentially the one Szeg\H{o} 
presented in \cite{Sz20,Sz21}. 

\smallskip
3. The same result is true for the symbol $d\mu = w(\theta)\f{d\theta}{2\pi} + d\mu_\s$, 
that is, the limit is independent of $d\mu_\s$. This extension was first proven by 
Verblunsky \cite{V36}. We will not prove this more general result here (\cite{OPUC} 
has four different proofs in Chapter~2) since it is peripheral to our main result. 

\smallskip
The first half of the theorem is a simple use of Jensen's inequality: 

\begin{proposition}\lb{P4.2} Let $w=e^L$ with $w,L\in L^1$. Then 
\begin{equation} \lb{4.2} 
\|\Phi_n\|_{w\f{d\theta}{2\pi}}^2 \geq \exp\biggl( \int L(\theta)\, 
\f{d\theta}{2\pi} \biggr)
\end{equation}
In particular, 
\begin{equation} \lb{4.3} 
D_n(w) \geq \exp\biggl( \int (n+1)L(\theta)\, \f{d\theta}{2\pi}\biggr)
\end{equation}
\end{proposition} 

\begin{proof} \eqref{4.3} follows from \eqref{4.2} and \eqref{2.3a}. To prove 
\eqref{4.2}, we write 
\begin{align} 
\|\Phi_n\|^2 = \|\Phi_n^*\|^2 &= \int \exp (2\log \abs{\Phi_n^*(e^{i\theta})} + 
L(e^{i\theta}))\, \f{d\theta}{2\pi} \notag \\
&\geq \exp\biggl( \int [2\log \abs{\Phi_n^*(e^{i\theta})} + L(e^{i\theta})]\biggr)\, 
\f{d\theta}{2\pi} \lb{4.4x} 
\end{align} 
by Jensen's inequality. 

By Theorem~\ref{T3.1}, $\log (\Phi_n^*(z))$ is analytic in $\bbD$, so 
\begin{align*} 
\int \log \abs{\Phi_n^* (e^{i\theta})}\, \f{d\theta}{2\pi} 
&= \Real \int \log (\Phi_n^*(e^{i\theta}))\, \f{d\theta}{2\pi} \\ 
&= \log \abs{\Phi_n^* (0)} =0 
\end{align*} 
since $\Phi_n$ monic implies $\Phi_n^*(0)=1$. 
\end{proof} 

The other half of the theorem depends on a variational principle noted by Szeg\H{o}: 

\begin{proposition}\lb{P4.3} We have for any $w\in L^1 (\partial\bbD, \f{d\theta}{2\pi})$, 
\begin{equation} \lb{4.4} 
\lim_{n\to\infty} \, [D_n (w)]^{1/n} = \inf \biggl\{\int \abs{f(e^{i\theta})}^2 
w(\theta)\, \f{d\theta}{2\pi} \biggm| f\in H^\infty (\bbD);\, f(0)=1 \biggr\}
\end{equation}
\end{proposition} 

{\it Remark.} $H^\infty (\bbD)$ is the Hardy space of bounded analytic functions on 
$\bbD$. By general principles \cite{Duren,Rudin}, for $\f{d\theta}{2\pi}$ a.e.~$e^{i\theta}
\in\partial\bbD$, $\lim_{r\uparrow 1} f(re^{i\theta})$ exists, and that is what we mean 
by $f(e^{i\theta})$ in \eqref{4.4}. 

\begin{proof} Since $\|\Phi_n^*\| = \|\Phi_n\|$, by \eqref{2.4}, 
\begin{align}  
\text{LHS of \eqref{4.4}} &= \lim_{n\to\infty}\, \|\Phi_n^*\|^2 \notag \\ 
&= \inf_n \, \|\Phi_n^*\|^2 \lb{4.5} 
\end{align} 
by \eqref{2.3b}. By the argument at the start of the proof of Theorem~\ref{T2.3}, 
\[
\Phi_n^* = \pi_n 1 
\]
where $\pi_n$ is the projection in the space of polynomials, $\calP_n$, of degree 
$n$ onto the orthogonal component of the span of $z,z^2,z^3, \dots, z^n$. This span 
is $\{P\in\calP_n \mid P(0) =0\}$, so by standard geometry, 
\begin{equation} \lb{4.6} 
\|\Phi_n^*\|^2 = \inf \{\|P\|^2 \mid P\in\calP_n;\, P(0)=1\}
\end{equation}
proving again that $\|\Phi_n^*\|^2$ is decreasing in $n$ and proving \eqref{4.4} if 
$H^\infty$ is replaced by the set of all polynomials. 

To complete the proof, we need only show that for any $f\in H^\infty$ with $f(0)=1$, 
there are polynomials $P_\ell(z)$ so that $P_\ell (0)=1$ and 
\begin{equation} \lb{4.7} 
\int\, \abs{P_\ell (e^{i\theta})}^2 w(\theta) \, \f{d\theta}{2\pi} \to 
\int\, \abs{f(e^{i\theta})}^2 w(\theta)\, \f{d\theta}{2\pi}
\end{equation}

If $f$ is analytic in a neighborhood of $\bbD$, the Taylor approximations converge 
uniformly on $\bar\bbD$, so \eqref{4.7} holds. For general $f$, by the dominated 
convergence theorem, 
\[
\lim_{r\uparrow 1} \, \int \, \abs{f(re^{i\theta})}^2 w(\theta)\, \f{d\theta}{2\pi} 
= \int \, \abs{f(e^{i\theta})}^2 w(\theta)\, \f{d\theta}{2\pi} 
\]
so, by a two-step approximation, we find $P_\ell$'s so \eqref{4.7} holds. 
\end{proof} 

\begin{proof}[Proof of Theorem~\ref{T4.1}] We will prove that for any $\veps >0$, there 
is an $f$ in $H^\infty$ with $f(0)=1$ and 
\begin{equation} \lb{4.7a} 
\int\, \abs{f(e^{i\theta})}^2 w(\theta) \, \f{d\theta}{2\pi} \leq \exp \biggl( \int 
\log (w(\theta)+\veps)\, \f{d\theta}{2\pi}\biggr)
\end{equation}
so taking $\veps\downarrow 0$ yields the opposite inequality to \eqref{4.3}. 

Define 
\begin{equation} \lb{4.8} 
g(z) = \exp \biggl( -\int \log (w(\theta)+\veps) \biggl( \f{e^{i\theta}+z}{e^{i\theta}-z} 
\biggr) \, \f{d\theta}{4\pi}\biggr)
\end{equation}
and $f(z) = g(z)/g(0)$, so $f(0)=1$. Moreover, $\abs{g(z)}\leq \veps^{-1/2}$ by the 
fact that the Poisson kernel 
\[
P_r (\theta, \varphi) = \Real \biggl( \f{e^{i\theta} + re^{i\varphi}}
{e^{i\theta} - re^{i\varphi}}\biggr) 
\]
is nonnegative and $\int \f{d\theta}{2\pi} P_r (\theta,\varphi)=1$. By standard maximal 
function arguments (see \cite{Rudin}), $\abs{g(e^{i\theta})}=\abs{w(\theta) +\veps}^{-1/2}$, 
so 
\[
\int \, \abs{f(e^{i\theta})}^2 w(\theta) \, \f{d\theta}{2\pi} \leq 
g(0)^{-2} = \text{RHS of \eqref{4.7a}} 
\]
\end{proof}

%%%%%%%%%%%%%%%%%%%%%%%%%%%%%%%%%%%%%%%%%%%
\section{The Szeg\H{o} Function} \lb{s5} 
%%%%%%%%%%%%%%%%%%%%%%%%%%%%%%%%%%%%%%%%%%%

When $w(\theta) =e^{L(\theta)}$ with $L\in L^1$, Szeg\H{o} \cite{Sz20,Sz21} introduced 
a natural function, $D(z)$ on $\bbD$, which will play a critical role in several places 
below: 
\begin{equation} \lb{5.1} 
D(z) =\exp\biggl( \int \biggl( \f{e^{i\theta} +z}{e^{i\theta}-z}\biggr) L(\theta)\,   
\f{d\theta}{4\pi}\biggr)
\end{equation}

Do not confuse $D_n$ and $D(z)$. Both symbols are standard, but the objects are very 
different.  

\begin{theorem}\lb{T5.1}  
\begin{SL} 
\item[{\rm{(a)}}] If \eqref{2.13} holds, then for $\abs{z}<1$, 
\begin{align} 
D(z) &= D(0) \exp\biggl( \, \sum_{k=1}^\infty \hat L_k z^k \biggr) \lb{5.2} \\ 
D(0) &= \biggl[ c_0 \biggl( w\, \f{d\theta}{2\pi}\biggr)\biggr]^{1/2} \,
\prod_{n=0}^\infty \, (1-\abs{\alpha_n}^2)^{1/2} \lb{5.3} 
\end{align} 
\item[{\rm{(b)}}] $D(z)$ lies in $H^2 (\bbD)$. 
\item[{\rm{(c)}}] $\lim_{r\uparrow 1} D(re^{i\theta})\equiv D(e^{i\theta})$ exist for 
a.e.~$\theta$ and 
\begin{equation} \lb{5.4} 
\abs{D(e^{i\theta})}^2 =w(\theta)
\end{equation} 
\item[{\rm{(d)}}] $D$ is nonvanishing on $\bbD$. 
\end{SL} 
\end{theorem} 

\begin{proof} (a) We get \eqref{5.2} and 
\begin{equation} \lb{5.5} 
D(0) =\exp (\tfrac12\, L_0)
\end{equation}
from \eqref{5.1} and 
\[
\f{e^{i\theta} + z}{e^{i\theta} -z} = 1+2\sum_{j=0}^\infty \, (e^{-i\theta} z)^n 
\]
uniformly in $e^{i\theta}\in\partial\bbD$ and $z\in\{\abs{z}<r\}$. \eqref{5.3} then 
follows from \eqref{2.11}, \eqref{4.1}, and \eqref{5.5}. 

\smallskip
(b) Let $D^{(M)}(z)$ be given by \eqref{5.1} with $L(\theta)$ replaced by $\min 
(L(\theta),M)$. Then $D^{(M)}\in H^\infty$ and $\abs{D^{(M)}(e^{i\theta})}^2 =
\min (w(\theta), e^{2M})$, 
\begin{align*} 
\sup_{0<r<1}\, \int_0^{2\pi} \, \abs{D^{(M)} (re^{i\theta})}^2 \, \f{d\theta}{2\pi} 
&= \int_0^{2\pi} \min (w(\theta), e^{2M})\, \f{d\theta}{2\pi} \\
&\leq \int_0^{2\pi} w(\theta)\, \f{d\theta}{2\pi} 
\end{align*} 
Thus, taking $M\to\infty$, we see $D\in H^2$. 

\smallskip
(c) is immediate from properties of boundary values of the Poisson integral in \eqref{5.1}. 

\smallskip
(d) is trivial from \eqref{5.1}.  
\end{proof} 

The following simple but powerful $L^2$ calculation goes back to Szeg\H{o} \cite{Sz20,Sz21} 
(it has a version when $d\mu_\s \neq 0$; see \cite{OPUC}): 
 
\begin{theorem} \lb{T5.2} Let $d\mu = w\f{d\theta}{2\pi}$ where $w=e^L$, $L\in L^1$. 
Then, as $n\to\infty$, 
\begin{alignat}{2} 
&\text{\rm{(i)}} \qquad && \int\, \abs{D\varphi_n^*(e^{i\theta})-1}^2 \, \f{d\theta}{2\pi}  
\to 0 \lb{5.6} \\
&\text{\rm{(ii)}} \qquad && \int\, \abs{\varphi_n^* (e^{i\theta})-D^{-1} (e^{i\theta})}^2\, 
d\mu \to 0 \lb{5.7} \\
&\text{\rm{(iii)}} \qquad && \varphi_n^*(z) \to D(z)^{-1} \lb{5.7a} 
\end{alignat} 
uniformly on compact subsets of $\bbD$. 
\end{theorem} 

\begin{proof} (i) By \eqref{5.4}, $\int \abs{D\varphi_n^*}^2\f{d\theta}{2\pi} = \int 
\abs{\varphi_n^*}^2\, d\mu =1$, so \eqref{5.6} is equivalent to 
\begin{equation} \lb{5.8} 
\Real \int \varphi_n^*(e^{i\theta}) D(e^{i\theta})\, \f{d\theta}{2\pi} \to 1 
\end{equation}
$D\varphi_n^*$ is in $H^2$, so the Cauchy integral formula applies, that is, 
\begin{align} 
\text{LHS of \eqref{5.8}} &= \varphi_n^*(0) D(0) \notag \\
&=\kappa_n \kappa_\infty^{-1} \lb{5.9} 
\end{align} 
since $\Phi_n^*(0)=1$, $\varphi_n =\kappa_n\Phi_n$ and \eqref{5.3} and \eqref{3.16a} 
imply $D(0)=\kappa_\infty^{-1}$. \eqref{5.9} implies \eqref{5.8}. 

\smallskip
(ii) is immediate from \eqref{5.6} and $d\mu =w\f{d\theta}{2\pi} =\abs{D}^2 
\f{d\theta}{2\pi}$ by \eqref{5.4}. 

\smallskip
(iii) \eqref{5.6} says $D\varphi_n^*\to 1$ in $\bbH^2$ and so, a fortiori, uniformly 
on compact subsets of $\bbD$. 
\end{proof} 

We need to extend this result to a neighborhood of $\bar\bbD$ when $L$ is real analytic. 
As a preliminary, we note: 

\begin{lemma}\lb{L5.3} Let $d\mu = w\f{d\theta}{2\pi}$ be a probability measure with 
$\log w\in L^1 (\f{d\theta}{2\pi})$. Then 
\begin{equation} \lb{5.10} 
\alpha_n = -\kappa_\infty \int \ol{\Phi_{n+1}(e^{i\theta})}\, D(e^{i\theta})^{-1} \, 
d\mu (\theta)
\end{equation} 
\end{lemma} 

\begin{proof} We will prove for $m\geq n+1$ that 
\begin{equation} \lb{5.11} 
\alpha_n = -\kappa_m \int \ol{\Phi_{n+1}(e^{i\theta})}\, \varphi_m^* (e^{i\theta})\, 
d\mu(\theta)
\end{equation} 
from which \eqref{5.10} follows from \eqref{5.7}. 

$\varphi_m^*$ is orthogonal to $\{z^\ell\}_{\ell=1}^m$, so if $P$ is any polynomial 
of degree at most $m$ with $P(0)=0$, 
\begin{equation} \lb{5.12} 
\int \ol{P(e^{i\theta})}\, \varphi_m^* (e^{i\theta})\, d\mu(\theta) =0
\end{equation}
\eqref{5.11} follows from \eqref{5.12} by taking 
\[
P(z) =\bar\alpha_n \varphi_m^* (z) + \kappa_m \Phi_{n+1}(z)
\]
which has $P(0)=\bar\alpha_n \kappa_m + \kappa_m (-\bar\alpha_n)=0$. 
\end{proof} 

The following is due to Nevai-Totik \cite{NT89}; it is needed in Section~\ref{s11}: 

\begin{theorem}\lb{T5.4} Suppose that $d\mu =e^L \f{d\theta}{2\pi}$ and $e^{i\theta} 
\mapsto L(\theta)$ is analytic in a neighborhood of $\partial\bbD$. Then $\varphi_n^* 
(z)\to D(z)^{-1}$ uniformly in a neighborhood of $\bar\bbD$. Moreover, the Verblunsky 
coefficients obey $\abs{\alpha_n} \leq C_2 e^{-An/2}$ for some $A>0$. 
\end{theorem} 

\begin{proof} Analyticity says $\abs{\hat L_k}\leq Ce^{-A\abs{k}}$ for some $A>0$. 
So, by \eqref{5.2}, $D(z)$ is analytic and nonvanishing in a disk of radius $e^A$. 
In particular, if 
\begin{equation} \lb{5.13} 
D(z)^{-1} =\sum_{j=0}^\infty d_{j,-1} z^j
\end{equation}
then 
\begin{equation} \lb{5.14} 
\abs{d_{j,-1}}\leq C_1 e^{-A\abs{j}/2}
\end{equation}

Plug \eqref{5.13} into \eqref{5.10} and note that  
\[
\int \ol{\Phi_{n+1} (e^{i\theta})}\, e^{ik\theta}\, d\mu(\theta) =0
\]
for $k=0,1,\dots, n$. Thus 
\begin{align} 
\abs{\alpha_n} &\leq \kappa_\infty \sum_{k=n+1}^\infty \, \abs{d_{k,-1}} \, 
\biggl| \int \ol{\Phi_{n+1}(e^{i\theta})}\, e^{ik\theta}\, d\mu \biggr| \notag \\
&\leq \kappa_\infty \|\Phi_{n+1}\| \sum_{k=n+1}^\infty\, \abs{d_{k,-1}} \notag \\
&\leq \kappa_\infty \sum_{k=n+1}^\infty \, \abs{d_{k,-1}} \lb{5.15} 
\end{align} 
since $\|\Phi_{n+1}\|\leq 1$. So, by \eqref{5.14}, 
\begin{equation} \lb{5.16} 
\abs{\alpha_n} \leq C_2 e^{-A\abs{n}/2}
\end{equation}

By \eqref{2.9} and $\abs{\Phi_n^* (e^{i\theta})} = \abs{\Phi_n (e^{i\theta})}$, we 
have 
\begin{align} 
\sup_\theta\, \abs{\Phi_{n+1}^* (e^{i\theta})} 
&\leq (1+\abs{\alpha_n}) \sup_\theta \, \abs{\Phi_n^* (e^{i\theta})} \notag \\
&\leq \prod_{j=0}^n \, (1+\abs{\alpha_j}) \notag \\
&\leq \exp \biggl( \, \sum_{j=0}^\infty \, \abs{\alpha_j}\biggr) = 
C_4 <\infty \lb{5.17} 
\end{align} 
by iteration. $C_4<\infty$ follows from \eqref{5.16}. Since $\Phi_n^*$ is analytic, 
we get 
\begin{equation} \lb{5.18} 
\sup_{z\in\bar\bbD}\, \abs{\Phi_n^*(z)} \leq C_4
\end{equation}
and thus, since $\Phi_n^*(z) = z^n \, \ol{\Phi_n (1/\bar z)}$, we get 
\begin{equation} \lb{5.19} 
z\in\bbC\backslash\bbD \Rightarrow \abs{\Phi_n (z)} \leq C_4 \abs{z}^n
\end{equation}

Returning to \eqref{2.9}, 
\begin{align*} 
\sum_{n=0}^\infty\, \abs{\Phi_{n+1}^*(z) - \Phi_n^*(z)} 
&\leq \sum_{n=0}^\infty\, \abs{\alpha_n} \, \abs{\Phi_n(z)} \\ 
&\leq C_4 C_2 \sum_{n=0}^\infty\, \abs{ze^{-A/2}}^n 
\end{align*} 
showing $\Phi_n^*$, and so $\varphi_n^*$, converges uniformly in $\{z\mid \abs{z} 
\leq e^{A/4}\}$. Since the limit is $D^{-1}$ in $\bbD$, it is $D^{-1}$ in this 
larger disk. 
\end{proof} 

Finally, in terms of $D$, we want to rewrite the second term in the Szeg\H{o} 
asymptotic formula \eqref{1.4}: 

\begin{theorem}\lb{T5.5} Let $d\mu =e^{L(\theta)} \f{d\theta}{2\pi}$ with $L\in 
L^1$. Let $\hat L_k$ be given by \eqref{1.3}. Then 
\begin{equation} \lb{5.20} 
\sum_{k=1}^\infty k \abs{\hat L_k}^2 = \f{1}{\pi} \int_{\abs{z}\leq 1} \,  
\abs{D(z)}^{-2} \biggl| \f{\partial D}{\partial z}\biggr|^2 \, d^2 z
\end{equation}
where both sides can be infinite. 
\end{theorem} 

\begin{proof} \eqref{5.20} follows by taking $r\uparrow 1$ in 
\begin{equation} \lb{5.21} 
\sum_{k=1}^\infty k\, \abs{\hat L_k}^2 r^{2k} =\f{1}{\pi} \int_{\abs{z} \leq r} \, 
\abs{D(z)}^{-2} \biggl| \f{\partial D}{\partial z}\biggr|^2 \, d^2 z
\end{equation}
(by using monotone convergence). To prove \eqref{5.21}, note that by \eqref{5.2}, 
\[
\log \biggl[ \f{D(z)}{D(0)}\biggr] = \sum_{k=1}^\infty \hat L_k z^k 
\]
converges uniformly in $\abs{z}<R$ and that $\abs{D}^{-2}\abs{\f{\partial D}{\partial z}}^2 
= \abs{\f{\partial}{\partial z}\log D(z)}^2$. Thus \eqref{5.21} follows from 
\[
\f{1}{\pi} \int_{\abs{z}\leq r} \bar z^{k-1} z^{\ell-1} \, d^2z = k^{-1} \delta_{k\ell} 
r^{2k}
\]
and 
\[
\f{\partial}{\partial z}\, \log D(z) = \sum_{k=1}^\infty k \hat L_k z^{k-1}
\]
\end{proof}

%%%%%%%%%%%%%%%%%%%%%%%%%%%%%%%%%%%%%%%%%%%%%%%%%%%%%%%%%%%%% 
\section{Extending the Strong Szeg\H{o} Theorem, Part I} \lb{s6}
%%%%%%%%%%%%%%%%%%%%%%%%%%%%%%%%%%%%%%%%%%%%%%%%%%%%%%%%%%%%% 

With the Szeg\H{o} function and Fatou's lemma, we have the tools for the 
Golinskii-Ibragimov \cite{GoIb} half of the extension theorem: 

\begin{theorem}[\cite{GoIb}]\lb{T6.1} Suppose for any BS measure, $d\mu$, we 
know that 
\begin{equation} \lb{6.1} 
G(d\mu) = \exp \biggl( \, \sum_{k=1}^\infty k \abs{\hat L_k}^2 \biggr)
\end{equation}
Then for any measure $d\mu = e^L \f{d\theta}{2\pi}$ with $L\in L^1$, we have 
\begin{equation} \lb{6.2} 
G(d\mu) \geq \exp \biggl( \, \sum_{k=1}^\infty k \abs{\hat L_k}^2\biggr)
\end{equation}
\end{theorem} 

{\it Remark.} In particular, $G<\infty \Rightarrow \sum k\abs{\hat L_k}^2 <\infty$. 

\begin{proof} Let $d\mu^{(N)}$ be the BS approximations to $d\mu$. Let $D^{(N)}(z)$ be 
the $D$ function for $d\mu^{(N)}$ and let 
\begin{equation} \lb{6.3}
M^{(N)} = \f{1}{\pi} \int_{\abs{z}\leq 1} \, \abs{D^{(N)}(z)}^{-2} \, 
\biggl| \f{\partial D^{(N)}}{\partial z}\biggr|^2 \, d^2 z
\end{equation}
and similarly for $M$ and $D$. 

By \eqref{2.12} and \eqref{3.20}, 
\[
G(d\mu^{(N)}) = \prod_{j=0}^{N-1}\, (1-\abs{\alpha_j}^2)^{-j-1} 
\]
so $G(d\mu^{(N)})$ is monotone increasing to $G(d\mu)$, that is, 
\begin{equation} \lb{6.4}
G(d\mu) = \lim_{N\to\infty}\, G(d\mu^{(N)})
\end{equation}

On the other hand, since $\varphi_{N+j} (z,d\mu^{(N)}) =z^j \varphi_N (z, d\mu^{(N)})$, 
$\varphi_{N+j}^* (z, d\mu^{(N}) =\varphi_N^* (z,d\mu^{(N)})=\varphi_N^* (z,d\mu)$, so 
by \eqref{5.7a}, 
\begin{equation} \lb{6.5}
D^{(N)}(z) = \varphi_N^* (z,d\mu)^{-1}
\end{equation}
and thus, by \eqref{5.7a} again, 
\begin{equation} \lb{6.6}
D^{(N)}(z)\to D(z)
\end{equation} 
uniformly on compacts. By analyticity, the same thing is true for derivatives, so 
\[
\abs{D(z)}^{-2}\, \biggl| \f{\partial D}{\partial z}\biggr|^2 = \lim_{N\to\infty} 
\, \abs{D^{(N)}(z)}^{-2}\, \biggl| \f{\partial D^{(N)}}{\partial z}\biggr|^2
\]
Thus, by Fatou's lemma, 
\begin{equation} \lb{6.7}
M \leq \liminf M^{(N)} 
\end{equation}
By hypothesis, 
\begin{equation} \lb{6.8}
G(d\mu^{(N)}) =\exp (M^{(N)})
\end{equation}

\eqref{6.8}, \eqref{6.7}, and \eqref{6.4} imply that 
\[
G(d\mu) \geq \exp (M)
\]
which is \eqref{6.2}. 
\end{proof}

%%%%%%%%%%%%%%%%%%%%%%%%%%%%%%%%%%%%%%%%%%%%%%%% 
\section{The Coulomb Gas Representation} \lb{s7} 
%%%%%%%%%%%%%%%%%%%%%%%%%%%%%%%%%%%%%%%%%%%%%%%% 

In this section, we will provide an integral formula for $D_n$ that will be critical 
in the next section. This formula appeared in Szeg\H{o}'s first paper \cite{Sz15} on 
asymptotics of Toeplitz determinants. He used it there for a minor technical purpose 
--- and for us, too, it plays a relatively minor role. That said, it plays a central 
role in two proofs of the strong Szeg\H{o} theorem and several applications. 

While we will not pursue the Coulomb gas picture, if one uses $\abs{z_j -z_k} =\exp 
(\log \abs{z_j-z_k})$, the formula we give for $D_n$ says it is the partition function 
of a two-dimensional gas, and this point of view is the basis of Johansson's proof 
\cite{Jo88}. It also explains some interest in Toeplitz matrices in the physics 
literature \cite{Len64,Len72,FH69}. 

Using Weyl's relation that Haar measure restricted to the classes of $\bbU(n+1)$, 
the group of $n\times n$ unitary matrices, is essentially $(2\pi)^{-n-1}[(n+1)!]^{-1} 
\prod_{k<j} \abs{e^{i\theta_k} - e^{i\theta_j}}^2\, d\theta_0 \dots d\theta_n$ with 
$\{e^{i\theta_j}\}_{j=0}^n$ the eigenvalues of $U\in \bbU(n+1)$, one can use the 
formula below to rewrite $D_n$ as an integral over $\bbU(n+1)$. This is both the 
starting point of the Bump-Diaconis \cite{BD} proof and of the many applications of 
Toeplitz matrices in the theory of random matrices \cite{MehBk}. 

It is, of course, well-known that $\prod_{k<j} (z_k-z_j)$ is a Vandermonde 
determinant. Expanding two such products, we get 
\begin{equation} \lb{7.1}
\biggl| \, \prod_{0\leq k < j \leq n}\, (z_k -z_j)\biggr|^2 = 
\sum_{\pi, \sigma\in\Sigma_{n+1}} \, (-1)^\pi (-1)^\sigma 
\prod_{j=0}^n \bar z_j^{\pi(j)} z_j^{\sigma(j)}
\end{equation}
where $\Sigma_{n+1}$ is the permutations of $\{0,\dots, n\}$ to itself. 

Setting $z_j =e^{i\theta_j}$ and integrating $d\mu (\theta_0) \dots d\mu (\theta_n)$, 
we get that 
\begin{equation} \lb{7.2}
\int \abs{\pi (z_k -z_j)}^2 \, d\mu(\theta_0) \dots d\mu (\theta_n) = 
\sum_{\pi,\sigma\in\Sigma_{n+1}} \, (-1)^\pi (-1)^\sigma \prod_{j=0}^n 
T_{\pi (j) \sigma (j)}
\end{equation}
where $T_{k\ell}$ are the matrix elements of the $(n+1)\times (n+1)$ Toeplitz 
matrix. Now 
\[
\prod_{j=0}^n T_{\pi(j)\sigma(j)} = \prod_{j=0}^n T_{j(\sigma \pi^{-1})(j)} 
\]
and $(-1)^\pi (-1)^\sigma =(-1)^{\sigma\pi^{-1}}$. Thus, summing over $\sigma$ for $\pi$ 
fixed and then over $\pi$, we see the right side of \eqref{7.2} is $(n+1)! D_n (d\mu)$. 
Specializing to $d\mu =e^L \f{d\theta}{2\pi}$, we have proven 

\begin{theorem}[Coulomb Gas Representation for $D_n$]\lb{T7.1} Let $e^L$, $L\in L^1 
(\partial\bbD, \f{d\theta}{2\pi})$. Then with $z_k =e^{i\theta_k}$, 
\[
D_n \biggl( e^L\, \f{d\theta}{2\pi}\biggr) = [(n+1)!]^{-1} \int_{\partial\bbD^{n+1}}\, 
\biggl|\, \prod_{0\leq k < j \leq n}\, (z_k-z_j)\biggr|^2 
\mathrm{e}^{\sum_{j=0}^n L(\theta_j)} \prod_{j=0}^n \, \f{d\theta_j}{2\pi} 
\]
\end{theorem} 

%%%%%%%%%%%%%%%%%%%%%%%%%%%%%%%%%%%%%%%%%%%%%%%%%%%%%%%%%%%%%
\section{Extending the Strong Szeg\H{o} Theorem, Part II} \lb{s8}
%%%%%%%%%%%%%%%%%%%%%%%%%%%%%%%%%%%%%%%%%%%%%%%%%%%%%%%%%%%%%

The Coulomb representation and Fatou's lemma are the tools for Johansson's half of 
the extension theorem. A measure of the form $e^L \f{d\theta}{2\pi}$ where $L$ is a 
real Laurent polynomial (i.e., $\sum_{k=-n}^n \hat L_k e^{ik\theta}$ with $\hat L_{-k} 
=\bar{\hat L}_k$) is called a {\it GI measure} (after its early use in \cite{GoIb}). If 
$d\mu = e^L \f{d\theta}{2\pi}$ with $L\in L^1$, we define the {\it GI approximations}, 
$d\mu_{(N)}$, by 
\begin{equation} \lb{8.1}
d\mu_{(N)} = \exp\biggl( \, \sum_{\abs{k}\leq N}\, \hat L_k e^{ik\theta}\biggr) 
\, \f{d\theta}{2\pi} \equiv\exp (L_{(N)}(\theta)) \, \f{d\theta}{2\pi}
\end{equation}

\begin{theorem}[\cite{Jo88}]\lb{T8.1} Suppose for any GI measure, $d\mu$, we know that 
\begin{equation} \lb{8.2}
G(d\mu) = \exp \biggl( \, \sum_{k=1}^\infty k \abs{\hat L_k}^2\biggr)
\end{equation}
Then for any measure $d\mu = e^L \f{d\theta}{2\pi}$ with $L\in L^1$, we have 
\begin{equation} \lb{8.3}
G(d\mu) \leq \exp\biggl(\, \sum_{k=1}^\infty k \abs{\hat L_k}^2\biggr)
\end{equation}
\end{theorem} 

{\it Remark.} In particular, if $\sum_{k=1}^\infty k \abs{\hat L_k}^2 <\infty$, then 
$G(d\mu) <\infty$. 

\begin{proof} Without loss (since we can multiply $d\mu$ by a constant), suppose 
$\hat L_0 =0$. Given $d\mu =e^L \f{d\theta}{2\pi}$, let $d\mu_{(N)}$ be its GI 
approximations and let 
\begin{equation} \lb{8.4} 
C_{(N),m} = [(n+1)!]^{-1} \int_{\bbD^{m+1}}\, \prod_{0\leq j < k \leq m}\, \abs{z_j -z_k}^2 
\mathrm{e}^{\sum_{j=0}^m L_{(N)}(\theta_j)} \prod_{j=0}^m \,\f{d\theta_j}{2\pi}
\end{equation}
so, of course, $C_{(N),m}=D_m (d\mu_{(N)})$ by Theorem~\ref{T7.1}. Let $C_m$ be the 
integral with $L_{(N)}$ replaced by $L$. 

There is nothing to prove in \eqref{8.3} if $\sum_{k=1}^\infty k \abs{\hat L_k}^2 =
\infty$, so suppose the sum is finite. Thus, $L_{(N)}\to L$ in $L^2$, and so, we can 
find a subsequence $L_{(N_j)}\to L$ pointwise. By Fatou's lemma, for each $m$, 
\begin{equation} \lb{8.5} 
C_m \leq \liminf_{j\to\infty}\,   C_{(N_j),m} 
\end{equation}

By Theorem~\ref{T7.1}, \eqref{2.12a}, $\hat L_0=0$, and the assumption for GI measures, 
\begin{align*} 
C_{(N_j),m} &= G_m (d\mu_{(N_j)}) \\ 
&\leq G(d\mu_{(N_j)}) \\
&= \exp \biggl(\, \sum_{k=1}^{N_j} k \abs{\hat L_k}^2\biggr) \\
&\leq \exp \biggl( \, \sum_{k=1}^\infty k\abs{\hat L_k}^2\biggr) 
\end{align*} 
By \eqref{8.5} and $\hat L_0 =0$, 
\[
G_m(d\mu) = C_m \leq \exp \biggl(\, \sum_{k=1}^\infty k\abs{\hat L_k}^2\biggr)
\]
Taking $m\to\infty$ yields \eqref{8.3}. 
\end{proof} 

Theorems~\ref{T6.1} and \ref{T8.1} reduce the proof of the sharp form of the strong 
Szeg\H{o} theorem to proving the result for BS and GI measures. In both cases, $d\mu 
=e^L \f{d\theta}{2\pi}$ with $e^{i\theta}\mapsto L(\theta)$ analytic in a 
neighborhood of $\partial\bbD$. In the last three sections, we will prove the 
strong Szeg\H{o} theorem in that case.

%%%%%%%%%%%%%%%%%%%%%%%%%%%%%%%%%%%%%%
\section{The CD Formula} \lb{s9} 
%%%%%%%%%%%%%%%%%%%%%%%%%%%%%%%%%%%%%

If $u_1,u_2$ solve $-u''_j + Vu_j =\lambda_j u_j$ with $u_j(0)=0$, then $(\bar\lambda_1 
-\lambda_2) \bar u_1 u_2 = (\bar u_1 u'_2 - u_2 \bar u'_1)'$, so 
\begin{equation} \lb{9.1} 
(\bar\lambda_1 -\lambda_2) \int_0^a \ol{u_1 (x)}\, u_2 (x)\, dx = \bar u_1 (a) 
u'_2 (a) - \bar u'_1 (a) u_2 (a)
\end{equation}
The Christoffel-Darboux formula is just the analog of this Wronskian relation that many 
undergraduates learn! 

\begin{theorem}[CD Formula]\lb{T9.1} Let 
\begin{equation} \lb{9.2} 
K_n (z,\zeta) =\sum_{j=0}^n \, \ol{\varphi_j(\zeta)}\, \varphi_j (z) 
\end{equation}
Then 
\begin{align} 
K_n (z,\zeta) &= \f{\ol{\varphi_{n+1}^*(\zeta)}\, \varphi_{n+1}^*(z) - 
\ol{\varphi_{n+1}(\zeta)}\, \varphi_{n+1}(z)}{1-\bar\zeta z} \lb{9.3} \\
&= \f{\ol{\varphi_n^*(\zeta)}\, \varphi_n^*(\zeta) - z\bar\zeta 
\, \ol{\varphi_n(\zeta)}\, \varphi_n(z)}{1-\bar\zeta z} \lb{9.4}
\end{align} 
\end{theorem} 

\begin{proof} By using \eqref{3.17d}/\eqref{3.17e} and their conjugates to write 
$\varphi_{n+1}$ and $\varphi_{n+1}^*$ in terms of $\varphi_n$ and $\varphi_n^*$, 
we get 
\begin{equation} \lb{9.5} 
\ol{\varphi_{n+1}^*(\zeta)}\, \varphi_{n+1}^*(z) - \ol{\varphi_{n+1}(\zeta)}\, 
\varphi_{n+1}(z) = \ol{\varphi_n^*(\zeta)}\, \varphi_n^*(z) - z\bar\zeta\, 
\ol{\varphi_n(\zeta)}\, \varphi_n(z)  
\end{equation}
where we used $\rho_n^{-2} (1-\abs{\alpha_n}^2) =1$, and that the $\varphi^* \varphi$ 
cross terms cancel. 

First of all, \eqref{9.5} says that \eqref{9.3} is equivalent to \eqref{9.4}. In addition, 
if we let $Q_n (z,\zeta)$ be what we know is the common value of the right sides of 
\eqref{9.3} and \eqref{9.4}, then subtracting \eqref{9.3} for $n-1$ from \eqref{9.4} for 
$n$, we see that 
\begin{align*}
Q_n (z,\zeta) - Q_{n-1} (z,\zeta) 
&= (1-z\bar\zeta) \biggl[\, \f{\ol{\varphi_n (\zeta)}\, \varphi_n(z)}{1-z\bar\zeta}\biggr] \\
&= K_n (z,\zeta) - K_{n-1} (z,\zeta) 
\end{align*} 
so we need only prove $Q_0 (z,\zeta)=K_0 (z,\zeta)$. Since \eqref{9.4} for $n=0$ is $1$ 
and that is $K_0 (z,\zeta)$, \eqref{9.3} is proven. 
\end{proof} 

\begin{corollary}\lb{C9.2} For $e^{i\theta}\in\partial\bbD$, 
\begin{equation} \lb{9.6} 
K_n (e^{i\theta}, e^{i\theta}) \equiv \sum_{j=0}^n \, \abs{\varphi_j (e^{i\theta})}^2 = 
\left. -\f{\partial}{\partial r}\, \abs{\varphi_{n+1}^* (re^{i\theta})}^2 \right|_{r=1} 
+ (n+1) \abs{\varphi_{n+1}^* (e^{i\theta})}^2
\end{equation}
\end{corollary} 

\begin{proof} In \eqref{9.3}, take $z=\zeta =re^{i\theta}$ and take $r\uparrow 1$. 
If we note $\abs{\varphi_{n+1} (re^{i\theta})}^2 =r^{2n+2} \abs{\varphi_{n+1}^* (r^{-1} 
e^{i\theta})}^2$, we get two terms: the one from the $-(r^{2n+2}-1)/(1-r^2)$ gives 
the $(n+1)\abs{\varphi_{n+1}^* (e^{i\theta})}^2$ term in \eqref{9.6}, and the other 
terms give the derivative in \eqref{9.6}. 
\end{proof} 

The CD formula for OPUC is due to Szeg\H{o} \cite{Szb}. Our proof is taken from Golinskii 
\cite{Golppt}. 

\smallskip 
{\it Remark.} Integrating both sides of \eqref{9.6} yields 
\begin{equation} \lb{9.7} 
\int \left. \f{\partial}{\partial r} \, \abs{\varphi_{n+1}^* (re^{i\theta})}^2 
\right|_{r=1} \, d\mu(\theta) =1 
\end{equation}

%%%%%%%%%%%%%%%%%%%%%%%%%%%%%%%%%%%%%%%%%%%%%%%%%%%%%%%%%%%%%%%%%%%%%%%%%
\section{The Feynman-Hellman Theorem for Toeplitz Determinants} \lb{s10}
%%%%%%%%%%%%%%%%%%%%%%%%%%%%%%%%%%%%%%%%%%%%%%%%%%%%%%%%%%%%%%%%%%%%%%%%%

Here we want to begin with a simple but useful formula for $\f{\partial}{\partial\lambda} 
\|\Phi_n\|_{d\mu_\lambda}^2$, where the measure $d\mu_\lambda$ depends smoothly on 
$\lambda$. This appeared in \cite{OPUC}, strongly motivated by closely related 
ideas of Baik et al. \cite{BDMZ}. 

\begin{theorem}\lb{T10.1} Let $d\mu_\lambda$ be a family of positive measures on 
$\partial\bbD$ which are $C^1$ in $\lambda$ in the sense that $\lambda\mapsto 
c_n (d\mu_\lambda)$ is $C^1$ for each $n$. Then $\|\Phi_n (\dott, d\mu_\lambda)
\|_{\mu_\lambda}$ is $C^1$ and 
\begin{equation} \lb{10.1} 
\f{d}{d\lambda}\, \log \|\Phi_n\|_{\mu_\lambda}^2 =\int \abs{\varphi_n (e^{i\theta},  
d\mu_\lambda)}^2\, \f{d\mu_\lambda (\theta)}{d\lambda}  
\end{equation}
\end{theorem} 

\begin{proof} Since $\f{d}{d\lambda} \log \|\Phi_n\|^2 = (\f{d}{d\lambda} \|\Phi_n\|^2)/ 
\|\Phi_n\|^2$, it suffices to prove that 
\begin{equation} \lb{10.2} 
\f{d}{d\lambda}\, \|\Phi_n\|_{\mu_\lambda}^2 = \int \, \abs{\Phi_n (e^{i\theta}, 
d\mu_\lambda)}^2\, \f{d\mu_\lambda (\theta)}{d\lambda} 
\end{equation}
Since 
\begin{equation} \lb{10.3} 
\|\Phi_\lambda\|^2 = \int \, \ol{\Phi_n (e^{i\theta}, d\mu_\lambda)}\, \Phi_n 
(e^{i\theta}, d\mu_\lambda)\, d\mu_\lambda 
\end{equation}
its derivative is a sum of three terms, of which one term is the right side of \eqref{10.2} 
and the other two are 
\begin{equation} \lb{10.4} 
\int \, \ol{\f{\partial}{\partial\lambda}\, \Phi_n (e^{i\theta}, d\mu_\lambda)}\, 
\Phi_n (e^{i\theta}, d\mu_\lambda)\, d\mu_\lambda 
\end{equation}
and its conjugate. 

But, for all $\lambda$, $\Phi_n$ is monic, so $\f{\partial}{\partial\lambda} \Phi_n$ is a 
polynomial of degree at most $n-1$ and so orthogonal to $\Phi_n$. It follows that the 
term in \eqref{10.4} is zero. 
\end{proof} 

{\it Remark.} Those familiar with the Feynman-Hellman theorem \cite{Thir} will see a 
metaphysical link to this proof. 

\smallskip
We want to do two things with \eqref{10.1}. First, we sum several logs in $\log D_n 
(d\mu) =\sum_{j=0}^n \log \|\Phi_j\|^2$ (by \eqref{2.3a}), and second, we make an 
explicit choice of $d\mu_\lambda$: 

\begin{theorem}\lb{T10.2} Let $w_t$ be a family of $L^1$ functions, $C^1$ in $t$ so 
$\log w_t$ is also $C^1$ in $t$, and $w_0=1$. Then 
\begin{equation} \lb{10.5} 
\begin{split}
\log D_n \biggl(w_1\, \f{d\theta}{2\pi}\biggr) &= (n+1) \log (\|\Phi_{n+1}\|_{t=1}^2) \\
&\quad -\int_0^1 \, dt \int \biggl[ \f{d}{dt}\, (\log w_t)\biggr] 
\biggl( \left. \f{\partial}{\partial r}\, \abs{\varphi_{n+1}^* (re^{i\theta};d\mu_t)}^2   
\right|_{r=1}\biggr)\, d\mu_t (\theta)
\end{split}
\end{equation}
\end{theorem} 

\begin{proof} We have that 
\[
\f{d}{dt}\, w_t \f{d\theta}{2\pi} = \f{w'_t}{w_t} \, d\mu_t = \f{d}{dt}\, 
(\log w_t)\, d\mu_t
\]
Using this, \eqref{2.3a}, \eqref{10.1}, and the definition \eqref{9.2} of $K$\!, 
we get 
\begin{equation} \lb{10.6} 
\f{d}{dt}\, \log D_n \biggl( w_t\, \f{d\theta}{2\pi} \biggr) = 
\int \biggl[ \f{d}{dt}\, (\log w_t) \biggr] K(e^{i\theta}, e^{i\theta}; d\mu_t)\, 
d\mu_t
\end{equation} 
Now use \eqref{9.6} to get two terms. If we integrate $dt$, the first term gives 
the integral in \eqref{10.5}. The second term can be integrated using \eqref{10.1} 
to give the first term in \eqref{10.5}.  
\end{proof}

%%%%%%%%%%%%%%%%%%%%%%%%%%%%%%%%%%%%%%%%%%
\section{Completion of the Proof}\lb{s11}
%%%%%%%%%%%%%%%%%%%%%%%%%%%%%%%%%%%%%%%%%% 

We complete the proof, following \cite{OPUC}, by showing: 

\begin{theorem} \lb{T11.1} Let $d\mu(\theta) =e^{L(\theta)}\f{d\theta}{2\pi}$ where 
$e^{i\theta}\mapsto L(\theta)$ is analytic in a neighborhood of $\partial\bbD$. Then 
\begin{equation} \lb{11.1} 
\lim_{n\to\infty} \, [\log D_n (d\mu) -(n+1) \hat L_0 ]= \sum_{k=1}^\infty 
k \abs{\hat L_k}^2 
\end{equation}
\end{theorem} 

\begin{proof} Without loss, we suppose $\int d\mu(\theta) =1$. We claim first that 
(under the analyticity assumption) 
\begin{equation} \lb{11.2} 
\lim_{n\to\infty}\, (n+1) [\log \|\Phi_{n+1}\|^2 -\hat L_0] =0
\end{equation}
For in terms of the Verblunsky coefficients, 
\[
\|\Phi_{n+1}\|^2 = \prod_{j=0}^n \, (1-\abs{\alpha_j}^2) 
\]
(by \eqref{2.8a}), while \eqref{4.1} and \eqref{2.11} say that 
\[
e^{\hat L_0} = \prod_{j=0}^\infty\, (1-\abs{\alpha_j}^2) 
\]
Thus 
\begin{align*}
\|\Phi_{n+1}\|^2 e^{-\hat L_0} &= \prod_{j=n+1}^\infty \, (1-\abs{\alpha_j}^2)^{-1} \\
&= 1 + O(e^{-An/2}) 
\end{align*} 
by \eqref{5.16}. This proves \eqref{11.2}. 

Thus, by \eqref{10.5}, \eqref{11.1} is equivalent to there being a choice $w_t$ with 
\begin{equation} \lb{11.3} 
-\int_0^1 dt \int \biggl( \f{d}{dt}\, \log w_t\biggr) \, \left. \f{\partial}{\partial r} 
\, \abs{\varphi_{n+1}^* (re^{i\theta};d\mu_t)}^2\right|_{r=1}\, d\mu_t (\theta) \to 
\sum_{k=1}^\infty k \abs{\hat L_k}^2 
\end{equation}

Define for $0\leq t\leq 1$, 
\begin{equation} \lb{11.4} 
c_t = \log \biggl[ \int e^{tL(\theta)}\, \f{d\theta}{2\pi}\biggr]
\end{equation}
and 
\begin{equation} \lb{11.5} 
w_t (\theta) = \exp (tL(\theta)-c_t) 
\end{equation} 
so $\int d\mu_t =1$ also. Thus, the Szeg\H{o} function for $\mu_t$ is given by 
\begin{align} 
\log D_t (z) &= t \log D(z) - \tfrac12\, c_t \lb{11.6}  \\
&= t \sum_{k=1}^\infty z^k \hat L_k + \tfrac12\, t \hat L_0 -\tfrac12\, c_t \lb{11.7} 
\end{align} 

It follows that $D_t(z)$ is analytic in a $t$-independent disk, and then the Verblunsky 
coefficients obey 
\[
\abs{\alpha_n (d\mu_t)} \leq Ce^{-An/2} 
\]
uniformly in $t\in [0,1]$, and then, from the proof of Theorem~\ref{T5.4}, that 
$\varphi_n^* (z,d\mu_t)\to D_t (z)^{-1}$ uniformly in $t\in [0,1]$ and $z$ in $\{z\mid 
\abs{z}\leq e^{A/4}\}$. Thus, the $\f{\partial}{\partial r}$ derivative in the integral 
on the left side of \eqref{11.3} converges to a $\f{\partial}{\partial r}$ derivative 
of $D_t$. 

We also write 
\begin{equation} \lb{11.7a}
\f{\partial}{\partial r}\, \abs{\varphi_{n+1}^*}^2  
= 2 \abs{\varphi_{n+1}^*}^2 \Real \, \f{\partial}{\partial r}\, \log (\varphi_{n+1}^*) 
\to -2 \abs{D_t}^{-2} \Real \biggl( t\, \f{\partial}{\partial r}\, \log D\biggr) 
\end{equation} 
and we can use $\abs{D_t}^{-2}\, d\mu_t = \f{d\theta}{2\pi}$. 

By \eqref{11.5},  
\begin{equation} \lb{11.8} 
\f{d}{dt}\, \log w_t = L(\theta) -\f{d}{dt}\, c_t 
\end{equation} 
Since the $\f{d}{dt} c_t$ is $\theta$-independent and by using \eqref{9.7}, we see the 
$\f{d}{dt}c_t$ can be integrated to $c_1 -c_0 =0$. Thus, in \eqref{11.3} we can replace 
$\f{d}{dt}\log w_t$ by $L(\theta)$. The result is that 
\begin{equation} \lb{11.9} 
\text{LHS of \eqref{11.3}} \to \int_0^1 dt \int 2t \, L(\theta)\, \Real \biggl[ 
\f{\partial}{\partial r} \, \log D \biggr] \, \f{d\theta}{2\pi}
\end{equation}
The only $t$-dependence is the $2t$ and $\int_0^1 2t\, dt =1$, so 
\begin{align*} 
\text{RHS of \eqref{11.9}} 
&= \int \biggl[\, \sum_{k=-\infty}^\infty \hat L_k e^{ik\theta}\biggr] 
\Real \biggl(\, \sum_{k=1}^\infty k\hat L_k e^{ik\theta}\biggr)\, \f{d\theta}{2\pi} \\
&= \int \biggl(\, \sum_{k=-\infty}^\infty \hat L_k e^{ik\theta}\biggr) 
\biggl[\tfrac12\, \sum_{k=-\infty}^\infty\, \abs{k} \hat L_k e^{ik\theta}\biggr]\, 
\f{d\theta}{2\pi} \\
&= \tfrac12\, \sum_{k=-\infty}^\infty \, \abs{k} \bar{\hat L}_k \hat L_k = 
\sum_1^\infty k\abs{\hat L_k}^2 
\end{align*} 
as claimed.
\end{proof}

\smallskip
 
%%%%%%%%%%%%%%%%%%%%%%%%%%%%%

\end{document}